\documentclass[12pt]{amsart}
\usepackage {amssymb,amsmath}
\usepackage{graphicx}
\usepackage{hyperref}
\makeatletter
\def\@cite#1#2{{\m@th\upshape\bfseries%
[{#1\if@tempswa{\m@th\upshape\mdseries, #2}\fi}]}}
\makeatother
\newtheorem{thm}{Theorem}[section]

\newtheorem{cor}[thm]{Corollary}
\newtheorem{prop}[thm]{Proposition}
\newtheorem{rem}[thm]{Remark}

\newtheorem{prob}[thm]{Problem}
\newtheorem{conj}[thm]{Conjecture}
\newcommand{\Prf}{\noindent\textbf{Proof.\ }}
\newcommand{\bx}{\strut\hfill$\blacksquare$\medbreak}
\newcommand{\upbx}{\vspace{-2.5\baselineskip}\newline\hbox{}\hfill$\blacksquare$\newline
\medbreak}

\newcommand{\ca}{\mathrm{C}^*}

\newenvironment{spmatrix}{\left(\begin{smallmatrix}}{\end{smallmatrix}\right)}

\newcommand{\bbC}{{\mathbb{C}}}

\newcommand{\bbR}{{\mathbb{R}}}

 \newcommand{\C}{{\mathcal{C}}}
 
 \newcommand{\E}{{\mathcal{E}}}

 \newcommand{\V}{{\mathcal{V}}}
 \newcommand{\W}{{\mathcal{W}}}

\newcommand{\upchi}{{\raise.35ex\hbox{$\chi$}}}
\newcommand{\qand}{\quad\text{and}\quad}

\newcommand{\qiff}{\quad\text{if and only if}\quad}
\newcommand{\rank}{\operatorname{rank}}

\newcommand{\Tr}{\operatorname{Tr}}

\def\bra#1{\langle #1}
\def\ket#1{|#1 \rangle}
\def\kb#1#2{|#1\rangle\!\langle #2 |}
\def\one{{\mathchoice{\rm 1\mskip-4mu l}{\rm 1\mskip-4mu l}{\rm 1\mskip-4.5mu l}{\rm
1\mskip-5mu l}}}

\begin{document}

\title[Higher-Rank Numerical Ranges]{Higher-Rank Numerical Ranges and Compression Problems}
\author[M.D. Choi, D.W. Kribs, K. {\.Z}yczkowski]{Man-Duen Choi$^1$,
David W. Kribs$^{2,3}$ and Karol {\.Z}yczkowski$^{4,5,6}$}
\thanks{2000 {\it Mathematics Subject Classification.} 15A60, 15A90, 47A12, 81P68.}
\thanks{{\it key words and phrases.} Hilbert space, Hermitian matrices, normal matrices,
higher-rank numerical range, compression-values, dilations,
quantum error correction.}

\address{$^1$Department of Mathematics, University of Toronto,
Toronto, Ontario, Canada M5S 2E4}
\address{$^2$Department of Mathematics and Statistics,
University of Guelph, Guelph, Ontario, Canada  N1G 2W1}
\address{$^3$Institute for Quantum Computing, University of Waterloo, Waterloo,
ON, CANADA N2L 3G1}
\address{$^4$Perimeter Institute for
Theoretical Physics, 31 Caroline St. North, Waterloo, ON, CANADA
N2L 2Y5}
\address{$^5$Institute of Physics, Jagiellonian University,
       ul. Reymonta 4, {30-059} Cracow, Poland}
\address{$^6$Center for Theoretical Physics, Polish Academy of Sciences,
   Al.~Lotnik{\'o}w 32/44, {02-668} Warsaw, Poland}

%
\begin{abstract}
We consider higher-rank versions of the standard numerical range
for matrices. A central motivation for this investigation comes
from quantum error correction. We develop the basic structure
theory for the higher-rank numerical ranges, and give a complete
description in the Hermitian case. We also consider associated
projection compression problems.
\end{abstract}
\maketitle

\section{Introduction}\label{S:intro}

In this paper we initiate the study of higher-rank versions of the
standard numerical range for matrices. A primary motivation for us
arises through the basic problem of error correction in quantum
computing. Specifically, the development of theoretical and
ultimately experimental techniques to overcome the errors
associated with quantum operations is central to continued
advances in quantum computing. As it turns out, the numerical
ranges $\Lambda_k(T)$, for $k > 1$, defined below are intimately
related to this problem of ``quantum error correction''. In the
paper \cite{CKZ05b} we give applications of the results from the
present paper to this problem.

Let $T$ be an $N\times N$ matrix with complex entries. For $k\geq
1$, define the {\it rank-$k$ numerical range of $T$} as the subset
$\Lambda_k(T)$ of the complex plane given by
\begin{eqnarray}\label{nrdefn}
\quad\,\,\, \Lambda_k(T) = \big\{ \lambda\in\bbC : PTP = \lambda P
\,\, {\rm for \,\, some \,\, rank-}k \,\,{\rm projection\,\,} P
\big\},
\end{eqnarray}
where we use the term ``projection'' to mean ``orthogonal
projection''. Observe that the numerical range of $T$ is obtained
as
\begin{eqnarray}
\Lambda_1(T) = W(T) = \{ \bra{T\psi}\ket{\psi} \, :\,
\ket{\psi}\in\bbC^N,\, ||\,\ket{\psi}\,|| = 1 \}.
\end{eqnarray}

In our analysis, it is desirable to explicitly find the scalars
$\lambda$ and the associated projections $P$ in
Eq.~(\ref{nrdefn}). Thus, this ``compression problem'' will be the
focus of this paper. A search of the substantive literature on
numerical ranges reveals connections with two lines of
investigation. The ``$k$-numerical range'' introduced by Halmos in
\cite{Hal67} is the set of all $\lambda$ that satisfy $\lambda =
\Tr (PTP)$ for some rank-$k$ projection $P$. Evidently this set
includes the set $k\, \Lambda_k(T)$, but the reverse inclusion
does not hold in general. The ``$k$th matrix numerical range''
studied by several authors consists of the set $W(k:T)$ of all
matrices $X^*TX$, where $X$ is an $N\times k$ matrix such that
$X^*X = I$. The higher-rank numerical ranges $\Lambda_k(T)$ can
alternatively be formulated as $\Lambda_k(T) = W(k:T) \cap \{
\lambda I_k : \lambda\in\bbC\}$. See
\cite{BoDu71,Arv72,LiTs91,Far93} as examples of other entrance
points into the literature on generalized notions of the numerical
range.

The rest of the paper is organized as follows. In the next section
($\S$~\ref{S:evalues}) we discuss the basic structure theory for
the sets $\Lambda_k(T)$. In particular, we derive an explicit
characterization of these sets for all Hermitian matrices. We
state a conjecture and an open problem in the case of normal
matrices. We discuss some lower dimensional cases in
$\S$~\ref{S:special}, and in the penultimate section
($\S$~\ref{S:projns}), we present a method for constructing the
associated compression projections that captures all possible
projections in the Hermitian case. In the context of quantum error
correction, projections that correspond to elements of
$\Lambda_k(T)$, for $k>1$, must be explicitly identified. For
instance, in the rank-two case, such projections correspond to
quantum bits of information, or ``qubits'', that can be corrected
after particular quantum operations act (see
$\S$~\ref{S:conclude}).

\section{Compression-Values}\label{S:evalues}

In this section we investigate the basic structure theory of the
sets $\Lambda_k(T)$. We shall refer to elements of $\Lambda_k(T)$
as ``compression-values'' for $T$, since $\lambda\in\Lambda_k(T)$
if and only if the $k\times k$ scalar matrix $\lambda I_k$ is the
compression of $T$ to a $k$-dimensional subspace. This means that
$T$ is unitarily equivalent to a $2\times 2$ block matrix of the
form
\begin{eqnarray}
T = \left( \begin{matrix} \lambda I_k & A \\ B & C
\end{matrix}\right),
\end{eqnarray}
where $A$ is a $k\times(N-k)$ matrix, $B$ is an $(N-k)\times k$
matrix, and $C$ is an $(N-k)\times(N-k)$ matrix. Equivalently, $T$
is a ``dilation'' of the scalar matrix $\lambda I_k$, or,
$T-\lambda I$ maps a $k$-dimensional subspace into its orthogonal
complement.

The following set inclusions may be readily verified:
\begin{eqnarray}\label{inclusions}
W(T) = \Lambda_1(T) \supseteq \Lambda_2(T) \supseteq \ldots
\supseteq \Lambda_N(T).
\end{eqnarray}
The following properties are also easily checked:
\begin{itemize}
\item[$(i)$] $ \Lambda_k(\alpha T + \beta I) = \alpha \Lambda_k(T)
+ \beta \quad \forall \alpha,\beta\in\bbC$. \item[$(ii)$]
$\Lambda_k(T^*) = \overline{\Lambda_k(T)}$. \item[$(iii)$]
$\Lambda_k(T) \subseteq \Lambda_k({\rm Re}\,T) + i\,
\Lambda_k({\rm Im}\,T)$. \item[$(iv)$] $\Lambda_k(T\oplus S)
\supseteq \Lambda_k(T) \cup \Lambda_k(S)$. \item[$(v)$]
$\Lambda_{k_1+k_2}(T\oplus S) \supseteq \Lambda_{k_1}(T) \cap
\Lambda_{k_2}(S)$.
\end{itemize}

The numerical range $W(T)= \Lambda_1(T)$ is a non-empty, compact
and convex subset of the plane that includes the spectrum of $T$.
If $T$ is normal, then $W(T)$ is the convex hull of the
eigenvalues for $T$. In particular, if $T$ is Hermitian, then
$W(T)$ is the closed interval of the real line determined by the
minimal and maximal eigenvalues of $T$. The higher-rank numerical
ranges can, of course, be empty.  But compactness still holds in
general. The proof of the following result is elementary, hence we
leave it to the interested reader.

\begin{prop}\label{cptcon}
Let $T$ be an $N\times N$ matrix and let $k\geq 1$. Then the
rank-$k$ numerical range $\Lambda_k(T)$ forms a compact set.
\end{prop}


Now we give a description of the higher-rank numerical range for
large values of $k$ relative to $N$.

\begin{prop}\label{largek}
Let $T$ be an $N\times N$ matrix and suppose that $2k > N$. Then
the rank-$k$ numerical range $\Lambda_k(T)$ is an empty set or a
singleton set. If $\Lambda_k(T) = \{\lambda_0\}$ is a singleton
set with $2k > N$, then $\lambda_0$ is an eigenvalue of geometric
multiplicity at least $2k-N$. In particular, $\Lambda_N(T)$ is
non-empty if and only if $T$ is a scalar matrix.
\end{prop}

\Prf Given $2k > N$, assume that $\Lambda_k(T)$ is non-empty, and
contains $\lambda_0\neq\lambda_1$. Let $P_0$, $P_1$ be the
corresponding rank $k$ projections. Then the projection $P =
P_0\wedge P_1$ onto the intersection of the ranges of these two
projections is non-zero and satisfies $\lambda_0 P = PTP =
\lambda_1 P$. This contradiction shows that $\Lambda_k(T)$ is a
singleton set when it is non-empty.

For the second claim, the equality $P(T-\lambda_0 I)P = 0$ implies
\begin{eqnarray}
T - \lambda_0 I = (I-P) (T - \lambda_0 I) + P (T -\lambda_0 I)
(I-P).
\end{eqnarray}
Hence, $ \rank(T-\lambda_0I) \leq 2 \rank (I-P) = 2N - 2k$, and
so,
\begin{eqnarray}
\ker(T-\lambda_0I) \geq N - (2N - 2k) = 2k - N.
\end{eqnarray}
\upbx


In the normal case the previous result yields more detailed
information for large values of $k$.

\begin{cor}
Let $T$ be an $N\times N$ normal matrix and suppose that $2k > N$.
Then the rank-$k$ numerical range $\Lambda_k(T)$ is an empty set
or a singleton set. In fact, the case $\Lambda_k(T) =
\{\lambda_0\}$ occurs if and only if there is a $(2N-2k)\times
(2N-2k)$ matrix $T_0$ such that $T$ is unitarily equivalent to $
\lambda_0 I_{2k-N} \oplus T_0, $ and $\lambda_0$ belongs to
$\Lambda_{N-k}(T_0)$.
\end{cor}

We now derive a general description of the rank-$k$ numerical
range in the Hermitian case for arbitrary $k$.

\begin{thm}\label{mainthm}
Let $A$ be an $N\times N$ Hermitian matrix with eigenvalues
(counting multiplicities) given by $a_1\leq a_2 \ldots \leq a_N$
and let $k\geq 1$ be a fixed integer with $1\leq k \leq N$. Then
the rank-$k$ numerical range $\Lambda_k(A)$ coincides with  $[a_k,
a_{N-k+1}]$ which is:
\begin{itemize}
\item[$(i)$] a non-degenerate closed interval if  $a_k <
a_{N-k+1}$ , \item[$(ii)$] a singleton set if $a_k =a_{N-k+1}$,
\item[$(iii)$] an empty set if $a_k > a_{N-k+1}$ .
\end{itemize}
Moreover, $\Lambda_k(A)$ coincides with the  intersection of the
numerical ranges $W(V^*A V)$, where $V$ runs through all
isometries $V: \bbC^{N-k+1}\rightarrow\bbC^N$.
\end{thm}

\Prf Let $\lambda\in\Lambda_k(A)$ and let $P_k$ be a rank-$k$
projection with $P_k A P_k = \lambda P_k$. If
$V:\bbC^{N-k+1}\rightarrow\bbC^N$ is an isometry, then the
subspace $P_k\bbC^N$ and the range space $VV^*(\bbC^{N})$ have
non-zero intersection. Thus, there exists a unit vector
$\ket{\psi} \in \bbC^N$ such that $\ket{\psi} = P_k\ket{\psi} = V
V^* \ket{\psi}$. Let $\ket{\psi'}$ be the unit vector in
$\bbC^{N-k+1}$ given by $\ket{\psi'} = V^* \ket{\psi}$. Then we
have
\begin{eqnarray}
\bra{V^*A V \psi'}\ket{\psi'} &=& \bra{A\psi}\ket{\psi} \\ &=&
\bra{P_kAP_k\psi}\ket{\psi} = \lambda \bra{P_k\psi}\ket{\psi} =
\lambda.
\end{eqnarray}
Hence we have shown that $\lambda$ belongs to $W(V^*AV)$. As $V:
\bbC^{N-k+1}\rightarrow\bbC^N$ was an arbitrary isometry, it
follows that $\Lambda_k(A)$ is contained in the intersection of
all such numerical ranges $W(V^*AV)$.

Next, let $\{\ket{i}: 1\leq i \leq N-k+1 \}$ be a fixed
orthonormal basis for $\bbC^{N-k+1}$ and let $\{\ket{\psi_i}\}$ be
an orthonormal basis for $\bbC^N$ of eigenvectors for $A$
corresponding to the eigenvalues $a_1,\ldots,a_N$. Consider two
linear isometries $V_1, V_2 :\bbC^{N-k+1}\rightarrow\bbC^N$
defined by $V_1(\ket{i}) = \ket{\psi_i}, V_2(\ket{i}) =
\ket{\psi_{N-i+1}}$.

Then $V_1^* A V_1$ and $V_2^* A V_2$ are operators on
$\bbC^{N-k+1}$ that are diagonal with respect to the basis
$\{\ket{i}\}$, and we have $W(V_1^* A V_1) = [a_1,a_{N-k+1}]$ and
$W(V_2^* A V_2) = [a_{k},a_N]$. It follows that
\begin{eqnarray}
\Lambda_k(A) \subseteq \bigcap_V W(V^*A V) &\subseteq& W(V_1^* A V_1)\bigcap W(V_2^* A V_2) \\
 &=& [a_k,a_{N-k+1}].
\end{eqnarray}

We complete the proof by showing $\Lambda_k(A)$ contains the set
$[a_k,a_{N-k+1}]$ when $a_k \leq a_{N-k+1}$. Suppose first that
$a_{N+1-k}> a_k$ (and so $2k \leq N$). Fix $\lambda$ in the
interval $[a_k,a_{N+1-k}]$. We shall directly construct a rank-$k$
projection $P_k$ such that $P_k A P_k = \lambda P_k$. Consider the
set of $k$ pairs $\{a_{k+1-j},a_{N-k+j}\}$, $1 \leq j \leq k$. As
a notational convenience we shall write $\{b_j,b_j'\}$ for the
ordered pair $\{a_{k+1-j},a_{N-k+j}\}$, and so $b_j> b_j^\prime$.
(The following construction may be easily modified for any joint
partition of the sets $\{a_N, \ldots, a_{N-k+1}\}$ and
$\{a_k,\ldots ,a_1\}$ into ordered pairs.)

We may write $A$, up to unitary equivalence, as a direct sum
\begin{eqnarray}
A = \big(\oplus_j A_j \big) \oplus B,
\end{eqnarray}
where each $A_j$ is a diagonal $2\times 2$ matrix with spectrum
$\{b_j,b_j^\prime\}$, and $B$ is either vacuous, or is the
diagonal matrix with diagonal entries
$\{a_{k+1},\ldots,a_{N-k}\}$. As $\lambda$ satisfies,
\begin{eqnarray}
\lambda \in [a_k,a_{N-k+1}] \subseteq [b_j,b_j^\prime] = W(A_j)
\quad \forall\, 1\leq j \leq k,
\end{eqnarray}
we may find angles $\theta_j$ such that
\begin{eqnarray}
\lambda = b_j\cos^2\theta_j + b_j' \sin^2\theta_j \quad \forall\,
1\leq j \leq k.
\end{eqnarray}
Now define an orthonormal set of $k$ vectors by
\begin{eqnarray}
\ket{\phi_j} = \cos\theta_j \ket{\psi_{N-k+j}} + \sin\theta_j
\ket{\psi_{k-j+1}}  \quad \forall\, 1\leq j \leq k,
\end{eqnarray}
and the rank-$k$ projection $P_k$ onto the subspace spanned by
these vectors;
\[
P = \kb{\phi_1}{\phi_1} + \kb{\phi_2}{\phi_2} + \ldots +
\kb{\phi_k}{\phi_k}.
\]

It follows that $P_k A P_k = \lambda P_k$. Indeed, observe that
for $1\leq j \leq k$ we have
\begin{eqnarray}
\bra{A\phi_1}\ket{\phi_j} &=& \cos\theta_1
\bra{A\psi_{N}}\ket{\phi_j} + \sin\theta_1\bra{A\psi_1}\ket{\phi_j} \\
&=& a_{N} \cos\theta_1 \bra{\psi_{N}}\ket{\phi_j} +
a_1\sin\theta_1\bra{\psi_1}\ket{\phi_j} \\
&=& b_1\cos\theta_1\cos\theta_j \delta_{N,N-1+j} +
b_1'\sin\theta_1\sin\theta_j\delta_{j,1}  \\ &=&
\lambda\delta_{j,1}.
\end{eqnarray}
Similarly, $\bra{A\phi_i}\ket{\phi_j}=\lambda\delta_{ij}$ for
$1\leq i,j\leq k$.

The remaining case is characterized by the constraint $\lambda:=
a_k = a_{N-k+1}$. If, in addition, $a_{N-k+2}>a_{k-1}$, then we
may split the sets $\{a_N,\ldots,a_{N-k+2}\}$ and
$\{a_{k-2},\ldots,a_1\}$ into pairs as above, and similarly define
$k-1$ vectors $\ket{\phi_1},\ldots,\ket{\phi_{k-1}}$. As the final
vector we can take $\ket{\phi_k}:=\ket{\psi_k}$, and define $P_k =
\sum_{j=1}^{k}\kb{\phi_j}{\phi_j}$. If $a_{N-k+2}=a_{k-1}$, but
$a_{N-k+3}>a_{k-2}$, then we will use $\ket{\psi_k}$ and
$\ket{\psi_{k-1}}$ as two of the vectors. This process may be
continued, if required, to account for degeneracies in the
spectrum of $A$ around the eigenvalue $a_k$, and construct a
rank-$k$ projection which yields $\lambda\in\Lambda_k(A)$. The
result now follows. \bx

\begin{figure} [htbp]
      \begin{center} \
    \includegraphics[width=12cm,angle=0]{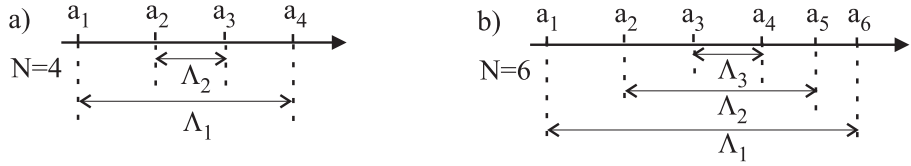}
\caption{The numerical range $\Lambda_k(A)$
   for a non-degenerate Hermitian operator $A$
   of size a) $N=4$, and b) $N=6$,
   with spectrum $\{ a_i \}$.}
\label{fig1}
\end{center}
    \end{figure}

For each real number $r\in \bbR$, we write $\lceil r \rceil$ for
the smallest integer $n$ satisfying $n\geq r$. From
Theorem~\ref{mainthm}, we see that if $k \leq \lceil N/2 \rceil$
(equivalently $2k-1\leq N$), then $\Lambda_k(A)$ is non-empty for
each $N\times N$ Hermitian matrix $A$. The following is an
analogous result for a general non-Hermitian matrix. (This result
can also be derived from Theorem~3.3 of \cite{LiTs91}.)

\begin{cor}\label{generalcor}
Let $T$ be an $N\times N$ complex matrix. Let $k$ be a positive
integer satisfying $k\leq \lceil N/4 \rceil$ (equivalently $4k-3
\leq N$). Then $\Lambda_k(T)$ is non-empty.
\end{cor}

\Prf Write $T = A + i B$ with $A=A^*$ and $B=B^*$. Let $b =
b_{2k-1}$ be the $(2k-1)$th smallest eigenvalue of $B$. By
Theorem~\ref{mainthm}, $b \in \Lambda_{2k-1}(B)$; and so there is
a projection $P$ of rank $2k-1$ such that $P (B-bI) P =0$.
Consider the $(2k-1)\times (2k-1)$ Hermitian matrix $A_0$ given by
the restriction of the compression $PAP$ to the range of $P$. It
follows from another application of Theorem~\ref{mainthm} that
$\Lambda_k(A_0)$ is a singleton set $\{a\}$, where $a$ is the
$k$th smallest eigenvalue of $A_0$. Hence there exists a
projection $Q \leq P$ such that $\rank Q = k$ and $QAQ = aQ$.
Thus, $QTQ = QAQ + i QBQ = (a+ib)Q$ and $\Lambda_k(T)$ is
non-empty.     \bx

The construction of projections that is described in the proof of
the previous theorem will be further fleshed out in subsequent
sections. It is perhaps appropriate to emphasize the most
important non-trivial case of this result. Specifically, when
$2k\leq N$ and the spectrum of $A$ is non-degenerate,
Theorem~\ref{mainthm} shows that the rank-$k$ numerical range is
the interval $\Lambda_k(A) = [a_{k},a_{N-k+1}]$
--- see Fig.~\ref{fig1}, which shows generalized numerical ranges
for $N=4$ and $N=6$. Also note that as an immediate consequence of
Theorem~\ref{mainthm}, it follows that the sets $\Lambda_k(A)$ are
convex for all $k\geq 1$ and Hermitian $A$.

We finish this section by discussing the case of normal matrices.
First note that property $(iii)$ above and Theorem~\ref{mainthm}
give a crude containment result for $\Lambda_k(T)$ for arbitrary
$T$. Indeed, $\Lambda_k(T)$ is a subset of the rectangular region
in the complex plane $\{\alpha + i \beta : \alpha\in\Lambda_k({\rm
Re}(T)),\,\beta\in\Lambda_k({\rm Im}(T))\}$. In general we can
obtain a more refined containment in the normal case.

\begin{thm}\label{normal}
Let $T$ be an $N\times N$ normal matrix and let $k\geq 1$ be a
fixed positive integer. Then
\begin{eqnarray}\label{normalid}
\Lambda_k(T) \subseteq \cap_\Gamma ({\rm co}\, \Gamma),
\end{eqnarray}
where $\Gamma$ runs through all $(N+1-k)$-point subsets (counting
multiplicities) of the spectrum of $T$.
\end{thm}

\Prf The relevant parts of the proof of Theorem~\ref{mainthm} can
be easily extended to the normal case to verify the inclusion of
Eq.~(\ref{normalid}).    \bx

\begin{rem}
{\rm Observe that Theorem~\ref{mainthm} shows the converse
inclusion of Eq.~(\ref{normalid}) holds in the Hermitian case. We
believe this inclusion holds more generally, at least in the
normal case, and we plan to undertake this investigation
elsewhere.}
\end{rem}

\begin{conj} {\rm If $T$ is an $N \times N$ normal matrix, then $\Lambda_k(T)$ coincides with the
intersection of the convex hulls ${\rm co}\,\Gamma$, where
$\Gamma$ is an $(N+1-k)$-point subset (counting multiplicities) of
the spectrum of $T$.}
\end{conj}

{\noindent}Verification of this conjecture would, of course,
automatically imply that $\Lambda_k(T)$ is convex, whenever this
set is non-empty and $T$ is normal. We state the general case as
an open problem.

\begin{prob}
{\rm Is $\Lambda_k(T)$ a convex set whenever it is non-empty?}
\end{prob}

\begin{figure} [htbp]
      \begin{center} \
    \includegraphics[width=11.0cm,angle=0]{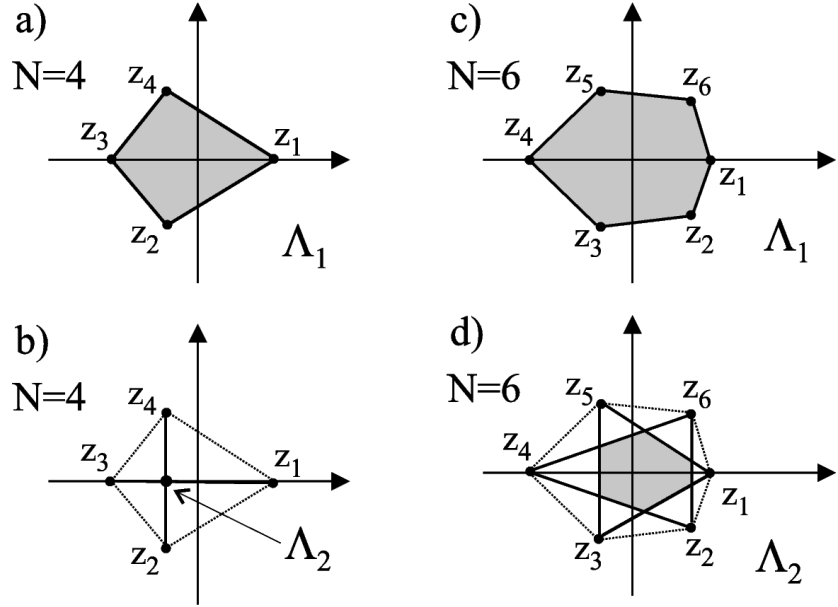}
\caption{Examples of the numerical range $\Lambda_k(T)$
   for a normal operator $T$
   of size a,b) $N=4$, and c,d) $N=6$,
   with non--degenerate complex eigenvalues $\{z_i\}$. The set $\Lambda_2$ is contained
   in the subset depicted (and equal at least in the case $N=4$).}
\label{fig2}
\end{center}
    \end{figure}

As a consequence of   Theorem~\ref{normal} and the proof of
Theorem~\ref{mainthm}, the conjecture can be seen to hold for
$N=1,2,3,4$ and all values of $k$ in each of these cases. Indeed,
Theorem~\ref{normal} shows that $\Lambda_k(T)$ is contained in the
desired set, and the construction of projections in the proof of
Theorem~\ref{mainthm} may be adapted for $N\leq 4$. In each of
these cases, $\Lambda_k(T)$ is either the empty set, a singleton
set, or an interval, and hence can never have interior. See Fig.
\ref{fig2} for an illustration of some of the non-interval cases
for $N=4$. (We note that the $N=4$ unitary case is explicitly
worked out in \cite{CKZ05b}.)

The first open case is that of $N=5$ and $k=2$. The cyclic 5-shift
is a good test example. This is the unitary
$U:\bbC^5\rightarrow\bbC^5$ defined on an orthonormal basis
$\{\ket{\xi_1},\ldots,\ket{\xi_5}\}$ by $U\ket{\xi_j} =
\ket{\xi_{j+1({\rm mod}\,5)}}$. The spectrum of $U$ is given by
$z_n = \exp\{i\frac{2\pi  n}{5}\}$, for $n=0,1,2,3,4$. Thus,
$\Lambda_2(U)$ is a subset of the pentagon shaped region depicted
in Fig.~3. The arguments of $(b)\Rightarrow(a)$ in
Theorem~\ref{mainthm} may be used to show that $\Lambda_2(U)$
contains the border points of this region, and also contains the
centre $\lambda =0$. The problem is to determine if the rest of
the interior points are included.

%
\begin{figure} [htbp]
       \begin{center} \
     \includegraphics[width=12.5cm,angle=0]{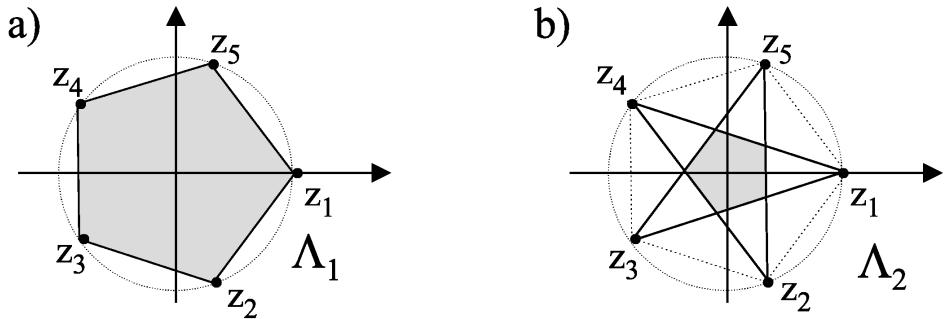}
\caption{Numerical ranges $\Lambda_k(U)$
    for the cyclic 5-shift,
  with spectrum consisting of the fifth roots of unity, $z_n$;
    a) $k=1$ and b) $k=2$, $\Lambda_2(U)$ is contained in this set.}
\label{fig3}
\end{center}
     \end{figure}

\section{Eigenvalue-Pairing Construction}\label{S:special}

The method presented in the next section shows how all the
higher-rank projections may be obtained through a generalization
of the ``eigenvalue pairing'' approach used in the proof of
Theorem~\ref{mainthm}. For illustration purposes, in this section
we further discuss the pairing approach in some lower dimensional
cases.

First let us recall the $k=1$ case as motivation for what follows
below. If $a_N\geq \ldots \geq a_1$ are the eigenvalues of $A=A^*$
as above, then the numerical range of $A$ is given by
$\Lambda_1(A) = [a_1,a_N]$. (Assume $A$ is non-scalar, so this is
truly an interval.) Let $\ket{\psi_j}$ be a choice of eigenvector
for each $a_j$. We may write a given $\lambda\in\Lambda_1(A)$ as a
linear combination $\lambda = a_1 c_1^2 + a_2 c_2^2 + \ldots + a_N
c_N^2$, where $\vec{c} = (c_i)$ are real scalars belonging to the
$(N-1)$-dimensional simplex ($\sum_{i=1}^N c_i^2 =1$). Then the
typical unit vector $\ket{\phi}$ that satisfies $\lambda =
\bra{A\phi}\ket{\phi}$, and the corresponding rank one projection,
are given by
\begin{eqnarray}
\ket{\phi} = \sum_{j=1}^N e^{i\theta_j} c_j\, \ket{\psi_j}
\quad\qand\quad P = \kb{\phi}{\phi}.
\end{eqnarray}
Observe that there are $N$ variables $c_j$, with two constraints,
and this gives $N-2$ free parameters. There is also an additional
$N$ free phases from the choices $\theta_j$, for $1\leq j \leq N$.
(Although it is $N-1$ free phases up to a global phase allowed in
the definition of $P$.)

For instance, if $(N,k) = (2,1)$, then $\Lambda_1(A) = [a_1,a_2]$
where $a_j$ is an eigenvalue for the eigenstate $\ket{\psi_j}$. In
this case, $\lambda\in [a_1,a_2]$ may be written as $\lambda = a_1
\cos^2 \beta_1 + a_2 \sin^2 \beta_1$, and the angle $\beta_1$ may
be computed via the equation
\begin{eqnarray}\label{2,1param}
\cos^2\beta_1 = \frac{\lambda - a_2}{a_1-a_2}.
\end{eqnarray}
The corresponding projection $P_1$ is obtained as a ``coherent
mixture'' of both eigenstates;
\begin{eqnarray}\label{2,1}
\ket{\phi_1} = e^{i\theta_1} \cos \beta_1 \ket{\psi_1} +
e^{i\theta_2} \sin \beta_1 \ket{\psi_2}, \quad P_1 =
\kb{\phi_1}{\phi_1}.
\end{eqnarray}

Next consider the case $(N,k)=(3,1)$. Let $\lambda$ belong to
$\Lambda_1(A) = [a_1,a_3]$. In this case, $\lambda$ may be
obtained via the equation
\begin{eqnarray}
\lambda = a_1 \cos^2\beta_1 +  a_2 \sin^2\beta_1 \cos^2\beta_2 +
a_3 \sin^2\beta_1 \sin^2\beta_2 .
\end{eqnarray}
In this case, $\beta_2 = \beta_2 (\lambda,\beta_1)$ depends on
both $\lambda$ and $\beta_1$, and hence there is a one parameter
family of solutions determined by $\beta_1$. The projection is
given by $P_1 = \kb{\phi_1}{\phi_1}$ where
\begin{eqnarray*}
\ket{\phi_1} = e^{i\theta_1} \cos\beta_1 \ket{\psi_1} +
e^{i\theta_2} \sin\beta_1\cos\beta_2 \ket{\psi_2} + e^{i\theta_3}
\sin\beta_1 \sin\beta_2 \ket{\psi_3},
\end{eqnarray*}
and we have three free phases $\{\theta_1,\theta_2,\theta_3\}$
(two phases up to a global phase). In the case that $\lambda =
a_2$, we may also use the solution Eq.~(\ref{2,1}) to find a
vector $\ket{\phi_{13}}$ as a mixture of $\ket{\psi_1}$ and
$\ket{\psi_3}$, and then mix it with $\ket{\psi_2}$ to obtain
\[
\ket{\phi_1} = \cos\beta_2\ket{\psi_2} + \sin\beta_2
\ket{\phi_{13}}.
\]

Let us turn now to higher-rank projections obtained from the
eigenvalue pairing approach in the case $N=4$. The case of
interest when $N=4$ is $(N,k)=(4,2)$. The challenge occurs when
$\Lambda_2(A) = [a_2,a_3]$ is a true interval. If we are given
$\lambda\in\Lambda_2(A)$, we can consider all pairs
$\{a_i,a_{i'}\}$ that contain $\lambda$. Here there are two
possibilities:
\begin{itemize}
\item[$(i)$] $\{a_4,a_2\}$, $\{a_3,a_1\}$, \item[$(ii)$]
$\{a_4,a_1\}$, $\{a_3,a_2\}$.
\end{itemize}
Of course, in the case of arbitrary $N$, there will be many more
possible pairings. Now we solve the $(2,1)$ problem for each of
the pairs separately. For instance, in the case of $(i)$, we solve
for $\beta_1$ and $\beta_2$ in the equations,
\begin{eqnarray}
\left\{ \begin{array}{rcl} \lambda &=& a_1\cos^2\beta_1 + a_3
\sin^2\beta_1 \\ \lambda &=& a_2\cos^2\beta_2 + a_4 \sin^2\beta_2,
\end{array}\right.
\end{eqnarray}
and so,
\begin{eqnarray}
\cos^2\beta_1 = \frac{\lambda-a_3}{a_1-a_3} \qand \cos^2\beta_1 =
\frac{\lambda-a_4}{a_2-a_4}.
\end{eqnarray}

We then define coherent combinations of eigenstates grouped in
pairs,
\begin{eqnarray}
\left\{ \begin{array}{rcl} \ket{\phi_1} &=&
e^{i\theta_1}\cos\beta_1\ket{\psi_1} + e^{i\theta_3}\sin\beta_1\ket{\psi_3} \\
\ket{\phi_2} &=& e^{i\theta_2}\cos\beta_2\ket{\psi_2} +
e^{i\theta_4}\sin\beta_2\ket{\psi_4}.
\end{array}\right.
\end{eqnarray}
Then write $P_2 = \kb{\phi_1}{\phi_1} +\kb{\phi_2}{\phi_2}$, and
it follows that $PAP=\lambda P$.

As noted above, this problem is equivalent to finding a unitary
matrix $U$ such that the matrix $A' = U A U^*$ includes a $2\times
2$ block given by the scalar matrix $\lambda \one_2$. In the case
$A = {\rm diag}\,(a_1,a_2,a_3,a_4)$, observe that one choice for
such a unitary is given by $U = O D$, where $D = {\rm
diag}\,(e^{i\theta_1},e^{i\theta_2},e^{i\theta_3},e^{i\theta_4})$
and $O$ is the orthogonal matrix given by
\begin{eqnarray}
O = \begin{spmatrix} \cos\beta_1 & 0 & \sin\beta_1 & 0 \\ 0 & \cos\beta_2 & 0 & \sin\beta_2  \\
- \sin\beta_1 & 0 & \cos\beta_1 & 0 \\ 0 & -\sin\beta_2 & 0 &
\cos\beta_2
\end{spmatrix} .
\end{eqnarray}

\section{Higher-Rank Projections}\label{S:projns}

In this section we consider the problem of finding all possible
rank-$k$ projections $P$ associated with a compression-value
$\lambda \in \Lambda_k(T)$; i.e., to solve for the rank-$k$
projections $P$ such that $PTP=\lambda P$. We shall focus on the
Hermitian case $A=A^*$. Recall that in the proof of
Theorem~\ref{mainthm} and the discussion of the previous section
we explicitly constructed certain families of projections to show
that particular values of $\lambda$ belonged to $\Lambda_k(A)$.
However, what we would like is a method for constructing such
projections that captures all possibilities. Unlike the standard
eigenvalue and eigenspace problem, in the generic case of this
compression problem there will be infinitely many projections.
Indeed, even in the typical case for the numerical range $W(A) =
\Lambda_1(A)$ this is the case. But it is possible, and in fact
easy, to write down such a method for the $k=1$ case. There are of
course more complications for $k\geq 2$.

First let us note that, while the eigenvalue-pairing approach
constructs a diverse set of projections, it is not sufficient to
capture all projections associated with values of $\Lambda_k(A)$.
Indeed, even consider the $k=1$ case of a Hermitian matrix $A$
with spectrum $\{0,1,2\}$. Here, $\Lambda_1(A) = [0,2]$. Let
$\ket{\psi_i}$, $i=0,1,2$, be unit eigenvectors for the
corresponding eigenvalues. The eigenvalue-pairing approach for
$\lambda=1$ in this case yields the family of projections
$P=\kb{\psi}{\psi}$, where
\begin{eqnarray}
\ket{\psi} = \frac{1}{\sqrt{2}}(e^{i\theta_1}\ket{\psi_0} +
e^{i\theta_2}\ket{\psi_2}).
\end{eqnarray}
But the set of all projections $P=\kb{\psi}{\psi}$ such that
$\bra{A\psi}\ket{\psi} = 1$ is the larger set given by unit
vectors of the form
\begin{eqnarray}
\ket{\psi} =
e^{i\theta_1}c_0\ket{\psi_0} + c_1\ket{\psi_1} + e^{i\theta_2}c_0
\ket{\psi_2}.
\end{eqnarray}

For an arbitrary Hermitian matrix $A$, the rank one projections
$P=\kb{\psi}{\psi}$ associated with values
$\lambda\in\Lambda_1(A)$ may be computed in the following manner.
Let $a_1\leq \ldots \leq a_N$ be the eigenvalues for $A$, and let
$\ket{\psi_i}$, $1\leq i \leq N$, be a choice of corresponding
eigenvectors. Suppose we have a unit vector $\ket{\psi} =
\sum_{i=1}^N c_i \ket{\psi_i}$. Then a simple computation shows
that
\begin{eqnarray}
\lambda = \bra{A\psi}\ket{\psi} \qiff \lambda = \sum_{i=1}^N a_i
|c_i|^2.
\end{eqnarray}

This constructive condition characterizes the rank one projections
associated with elements of the numerical range. Notice that there
are infinitely many possibilities for such projections whenever
$\lambda$ is not an eigenvalue for $A$. There is a corresponding
characterization for arbitrary $k$, though it is not constructive
for $k\geq 2$. Instead, in what follows  we present a
constructive, algorithmic approach to find all higher-rank
projections associated with compression-values of $\Lambda_k(A)$
for $k\geq 2$.

Let $\lambda\in\Lambda_k(A)$. By using a translation, we may
assume that $\lambda =0$. Let $P_+$, $P_0$, and $P_-$ be the
projections onto the eigenspaces of $A$ for respectively, the
positive eigenvalues, the eigenvalue zero, and the negative
eigenvalues. First, we consider the case when there is no
degeneracy in the spectrum of $A$; that is, $\lambda = 0$ is not
an eigenvalue of $A$.  Next choose a $k$-dimensional subspace
$\V_+$ of $P_+\bbC^N$. Note that this is possible by
Theorem~\ref{mainthm}, and our assumptions $\lambda =0\in
\Lambda_k(A) = [a_k,a_{N-k+1}]$, and the non-degeneracy of the
spectrum. By the same reasoning, $P_-\bbC^N$ is at least
$k$-dimensional, and hence we may choose an isometry $U: \V_+
\rightarrow P_-\bbC^N$. Now we define a $k$-dimensional subspace
of $\bbC^N$;
\begin{eqnarray}\label{Vform}
\V = \big\{ v_+ + U v_+ : v_+ \in \V_+ \big\}.
\end{eqnarray}


Next define a $k$-dimensional subspace $\W = f(A) \V$ where
\begin{eqnarray}
f(x) = \left\{ \begin{array}{cl}  |x|^{-1/2} & \mbox{if $x\neq 0$}
\\ 1 & \mbox{if $x=0$} \end{array} \right..
\end{eqnarray}
Observe that $f(A)Af(A) = P_+ - P_-$, and hence $\forall
v_1,v_2\in\V$ we have
\begin{eqnarray}
\bra{A f(A) v_1}\ket{f(A) v_2}  &=& \bra{f(A) A f(A) v_1}\ket{
v_2} \\ &=& \bra{(P_+ - P_- ) v_1}\ket{v_2} \\ &=& \bra{P_+
v_1}\ket{P_+v_2} - \bra{P_- v_1}\ket{P_- v_2} = 0.
\end{eqnarray}
It follows that $P_\W$ is a rank-$k$ projection such that $P_\W A
P_\W = 0$.

Now we show that every rank-$k$ projection $P$ such that $PAP=0$,
can be written in the form $P=P_\W$ as above. Let $P$ be such a
projection, and let $\V$ be the $k$-dimensional subspace $\V =
f(A)^{-1} \W$. Then for all $v\in\V$ we have
\begin{eqnarray}
0&=& \bra{A f(A) v}\ket{f(A) v}  = \bra{f(A) A f(A) v}\ket{ v}
\\ &=& \bra{(P_+ - P_- ) v}\ket{v} = || P_+ v||^2 - ||P_- v||^2.
\end{eqnarray}
In particular, this implies that the map $U(P_+ v) \equiv P_- v$
determines a well defined isometry $U: \V_+ \rightarrow \V_-$,
where $\V_+ = P_+ \V$ and $\V_- = P_-\V$. Thus, $\V_+$ and $\V_-$
are both $k$-dimensional and $\V$ is of the form given in
Eq.~(\ref{Vform}), and hence $P=P_\W$ as claimed.

We have presented a constructive method to obtain projections
associated with the compression-values of $\Lambda_k(A)$, in the
case that there are no degeneracies in the spectrum of $A$. We
have also shown that every such projection arises in this manner.
Let us summarize the method.
\begin{itemize}
\item[$(i)$] Choose a $k$-dimensional subspace $\V_+$ of
$P_+\bbC^N$. \item[$(ii)$] Choose a linear isometry $U:\V_+
\rightarrow P_-\bbC^N$. \item[$(iii)$] Define the $k$-dimensional
subspace $\V = \{ v_+ + Uv_+ : v_+ \in \V_+\}$. \item[$(iv)$] Let
$\W = f(A) \V$. Then $\dim\W =k$ and $P_\W A P_\W = 0$.
\end{itemize}

If there are degeneracies in the spectrum of $A$, the above method
may be adjusted by including part of the subspace $P_0\bbC^N$ in
the subspaces $\V_+$ and $\V_-$ as follows: As above, we want to
construct all $k$-dimensional subspaces $\V$ of $\bbC^N$ such that
\begin{eqnarray}
\bra{(P_+ - P_-) v}\ket{v} = 0 \quad \forall\, v\in\V.
\end{eqnarray}
This can be accomplished since $0 \in \Lambda_k(A) =
[a_k,a_{N-k+1}]$, and so $k\leq \dim P_0 \bbC^N + \dim
P_{\pm}\bbC^N$. Consider all possible pairs of non-zero integers
$(k_1,k_2)$ with $k_1+k_2 = k$, $k_1\leq \dim P_0\bbC^N$, and
\begin{eqnarray}
k_2 \leq \min \{\dim P_-\bbC^N, \dim P_+ \bbC^N \}.
\end{eqnarray}
Choose $\V_0$ as any $k_1$-dimensional subspace of $P_0 \bbC^N$
and choose $\V_+$ as any $k_2$-dimensional subspace of $P_+
\bbC^N$. Let $X: \V_+ \rightarrow (\V_0)^\perp \cap P_0\bbC^N$ be
any operator and let $U : \V_+ \rightarrow P_- (\bbC^N)$ be any
isometry and define
\begin{eqnarray}
\V = \V_0 + \{ v + Xv + Uv : v\in\V_+ \}.
\end{eqnarray}
Then $\W = f(A) \V$ is a $k$-dimensional subspace of $\bbC^N$ with
the desired properties.

\section{Concluding Remark}\label{S:conclude}

We conclude by briefly discussing the mathematical  context of the
work \cite{CKZ05b}, which includes applications of the present
work to quantum computing. Every quantum operation $\E$ on a given
quantum system is determined operationally by a set of operators
$\{ A_i \}$ that act on the Hilbert space for the system via the
so-called operator-sum representation $\E(\cdot) = \sum_i A_i
(\cdot) A_i^*$. (See \cite{Kri05a} for a brief introduction to
some of the mathematical aspects of quantum computing.) In the
context of quantum error correction, the $A_i$ are often called
``error operators''. It is the effects of such operators that must
be mitigated for whenever there is a transfer of quantum
information determined by $\E$. There are numerous strategies that
have been, and are being, developed for this type of error
correction. We go into detail on this subject in \cite{CKZ05b},
but here we indicate how the mathematical conditions that
characterize correction in the fundamental protocol for quantum
error correction connects with the higher-rank numerical ranges.
In the ``standard model'' for quantum error correction
\cite{BDSW96a,KL97a}, codes are identified with subspaces of the
system Hilbert space, and ``correctability'' of a given code
subspace $\C$ in terms of an error model $\E$ is shown to be
equivalent to the existence of scalars $\Lambda =(\lambda_{ij})$
such that
\begin{eqnarray}
P_\C A_i^* A_j P_\C = \lambda_{ij} P_\C \quad \forall \, i,j.
\end{eqnarray}
Here $P_\C$ denotes the projection of the system space onto $\C$.
Thus, the problem of finding correctable codes for a given error
model $\E = \{A_i\}$ is equivalent to  finding the
compression-values inside the higher-rank numerical ranges
$\Lambda_k(A_i^* A_j)$, $\forall i,j$ and $\forall k> 1$, along
with the corresponding projections. As indicated in \cite{CKZ05b},
this problem may be reduced to a system of such problems for
Hermitian or normal operators.



{\noindent}{\it Acknowledgements.} We are grateful to the referees
for helpful comments. D.W.K. would like to thank Ruben
Martinez-Avendano for a stimulating conversation on numerical
ranges. M.D.C. and D.W.K. were partially supported by NSERC. K.
{\.Z}. acknowledges a partial support by the grant number
PBZ-MIN-008/P03/2003 of Polish Ministry of Science and Information
Technology.


%


\begin{thebibliography}{30}
\expandafter\ifx\csname
natexlab\endcsname\relax\def\natexlab#1{#1}\fi
\expandafter\ifx\csname bibnamefont\endcsname\relax
  \def\bibnamefont#1{#1}\fi
\expandafter\ifx\csname bibfnamefont\endcsname\relax
  \def\bibfnamefont#1{#1}\fi
\expandafter\ifx\csname citenamefont\endcsname\relax
  \def\citenamefont#1{#1}\fi
\expandafter\ifx\csname url\endcsname\relax
  \def\url#1{\texttt{#1}}\fi
\expandafter\ifx\csname
urlprefix\endcsname\relax\def\urlprefix{URL }\fi
\providecommand{\bibinfo}[2]{#2}
\providecommand{\eprint}[2][]{\url{#2}}


\bibitem{CKZ05b}
\bibinfo{author}{\bibfnamefont{M.~D.}~\bibnamefont{Choi}},
  \bibinfo{author}{\bibfnamefont{D.~W.}~\bibnamefont{Kribs}}, \bibnamefont{and}
  \bibinfo{author}{\bibfnamefont{K.}~\bibnamefont{{\.Z}yczkowski}},
   \bibinfo{journal}{{\it Quantum error correction
               and higher-rank numerical range},}
               \bibinfo{pages}{quant-ph/0511101}.


\bibitem{Hal67}
\bibinfo{author}{\bibfnamefont{P.}~\bibnamefont{Halmos}},
   \bibinfo{journal}{{\it A Hilbert space problem book},}
   \bibinfo{journal}{D. Van Nostrand Company, Ltd., Toronto,} (\bibinfo{year}{1967}).



\bibitem{BoDu71}
\bibinfo{author}{\bibfnamefont{F.~F.}~\bibnamefont{Bonsall}} \bibnamefont{and}
  \bibinfo{author}{\bibfnamefont{J.}~\bibnamefont{Duncan}},
   \bibinfo{journal}{{\it Numerical ranges of operators on normed spaces
   and of elements of normed algebras},}\bibinfo{journal}{London Mathematical
   Society Lecture Note Series, 2, Cambridge University Press, London-New York} (\bibinfo{year}{1971}).

\bibitem{Arv72}
\bibinfo{author}{\bibfnamefont{W.}~\bibnamefont{Arveson}},
   \bibinfo{journal}{{\it Subalgebras of $\ca$-algebras, II},}
   \bibinfo{journal}{Acta Math.}  \textbf{\bibinfo{volume}{128}},
   \bibinfo{pages}{271-308} (\bibinfo{year}{1972}).

\bibitem{LiTs91}
\bibinfo{author}{\bibfnamefont{C.-K.}~\bibnamefont{Li}} \bibnamefont{and}
  \bibinfo{author}{\bibfnamefont{N.-K.}~\bibnamefont{Tsing}},
   \bibinfo{journal}{{\it On the $k$th matrix numerical range},}
   \bibinfo{journal}{Linear and Multilinear Algebra}  \textbf{\bibinfo{volume}{28}},
   \bibinfo{pages}{229-239} (\bibinfo{year}{1991}).

\bibitem{Far93}
\bibinfo{author}{\bibfnamefont{D.~R.}~\bibnamefont{Farenick}},
   \bibinfo{journal}{{\it Matricial extensions of the numerical range: A brief survey},}
   \bibinfo{journal}{Linear and Multilinear Algebra}  \textbf{\bibinfo{volume}{34}},
   \bibinfo{pages}{197-211} (\bibinfo{year}{1993}).

\bibitem{Kri05a}
\bibinfo{author}{\bibfnamefont{D.~W.}~\bibnamefont{Kribs}},
   \bibinfo{journal}{{\it A quantum computing primer for operator theorists},}
   \bibinfo{journal}{Linear Algebra Appl.}  \textbf{\bibinfo{volume}{400}},
   \bibinfo{pages}{147-167} (\bibinfo{year}{2005}).

  \bibitem{BDSW96a}
\bibinfo{author}{\bibfnamefont{C.~H.} \bibnamefont{Bennett}},
  \bibinfo{author}{\bibfnamefont{D.~P.} \bibnamefont{DiVincenzo}},
  \bibinfo{author}{\bibfnamefont{J.~A.} \bibnamefont{Smolin}},
  \bibnamefont{and} \bibinfo{author}{\bibfnamefont{W.~K.}
  \bibnamefont{Wootters}},  \bibinfo{journal}{{\it Mixed-state entanglement and quantum error
  correction},}
   \bibinfo{journal}{Phys. Rev. A}
  \textbf{\bibinfo{volume}{54}}, \bibinfo{pages}{3824} (\bibinfo{year}{1996}).

\bibitem{KL97a}
\bibinfo{author}{\bibfnamefont{E.}~\bibnamefont{Knill}} \bibnamefont{and}
  \bibinfo{author}{\bibfnamefont{R.}~\bibnamefont{Laflamme}},
   \bibinfo{journal}{{\it Theory of quantum error-correcting codes},}
   \bibinfo{journal}{Phys. Rev. {A}} \textbf{\bibinfo{volume}{55}},
  \bibinfo{pages}{900} (\bibinfo{year}{1997}).




\end{thebibliography}
\end{document}